\documentclass{amsart}
\usepackage{amssymb,amsmath,latexsym,amsthm}
\usepackage[english]{babel}

\usepackage{amsmath,amssymb,amsbsy,amsfonts,amsthm,latexsym,
                       amsopn,amstext,amsxtra,euscript,amscd}
\usepackage{hyperref}
\usepackage{array}
\usepackage{ifthen}
\usepackage{url}

\newtheorem{thm}{Theorem}

\newtheorem{lem}[thm]{Lemma}

\def\NN{\mathbb N}

\begin{document}

\title[On the Exponential Diophantine Equation $(a^n-1)(b^n-1)=x^2$]{On the Exponential Diophantine Equation $(a^n-1)(b^n-1)=x^2$}

\author[Armand Noubissie]{Armand Noubissie}
\address{Institut de Math\'ematiques et de Sciences Physiques. Dangbo, B\'enin}
\email{armand.noubissie@imsp-uac.org}  

\author[Alain Togb\'e]{Alain Togb\'e}
\address{Department of Mathematics, Statistics and Computer Science,
Purdue University Northwest, 1401 S, U.S. 421, Westville IN 46391 USA}
\email{atogbe@pnw.edu}

\author[Zhongfeng Zhang]{Zhongfeng Zhang}
\address{School of Mathematics and Statistics, Zhaoqing University, Zhaoqing 526061, P. R. China}
\email{bee2357@163.com}



\date{\today}

\maketitle

\begin{abstract}
Let $a$ and $b$ be two distinct fixed positive integers such that $\min \{a,b\}>1.$  First, we  correct an oversight from \cite{X-Z}. Then, we show that the equation in the title with $b \equiv 3 \pmod 8$, $b$ prime and $a$ even has no solution in positive integers $n, x$. This generalizes a result of Szalay \cite{L}. 

\end{abstract}

\section{Introduction}\label{sec1}

Let $\mathbb{N}^+$ be the set of all positive integers. Let $a>1$ and  $b>1$ be different fixed integers. The exponential Diophantine equation 
\begin{equation}\label{eq1}
(a^n-1)(b^n-1)=x^2,~ x,n \in \mathbb{N}^+  
\end{equation} 
has been studied by many authors in the literature since 2000. First, Szalay \cite{L} studied equation \eqref{eq1} for $(a,b)=(2,3)$ and showed that this equation has no positive integer solutions. Our Theorem \ref{thm2} generalizes this result. He also proved that equation \eqref{eq1} has only the positive integer solution $(n,x)=(1,2)$, for $(a,b)=(2,5)$ and there are no solutions, for $(a,b)=(2,2^k)$ with $k\geq 2$ except for $n=3$ and $k=2$. Hajdu and Szalay \cite{H-L} proved that \eqref{eq1} has no solution for $(a,b)=(2,6)$ and for  $(a,b)=(a,a^k),$ there are no solutions with $k\geq 2$ and $kn> 2$ except for the three cases $(a,n,k)=(2,3,2), (3,1,5),(7,1,4)$ generalizing Theorem 3 of \cite{L}. This result was extended by Cohn \cite{C} to the case $a^k=b^l$ (see RESULT 1). Cohn also proved that there are no solutions to \eqref{eq1} when $4\mid n$, except for $\{a,b\}= \{13,239\}$ with $n=4.$ Luca and Walsh  \cite{W-L} proved that equation \eqref{eq1} has finitely positive solutions for fixed $(a,b)$ and completely solved the equation for almost all pairs $(a,b)$ in the range $1<a<b\leq 100.$ Since then, many authors studied equation \eqref{eq1} by introducing some special contraints on $a$ or $b$ (see, for example  \cite{I}, \cite{L-L}, \cite{X-Z}).

One of our main results is the following.

\begin{thm} \label{thm1}
Let $a,~b \in \mathbb{N}$ such that $a,~b>1.$ Suppose that one of the following properties is satisfied:
\begin{itemize}
\item $a \equiv 2 \pmod 3$ and $b \equiv 0 \pmod 3$;
\item $a \equiv 3 \pmod 4$ and $b \equiv 0 \pmod 2$;
\item $a \equiv 4 \pmod 5$ and $b \equiv 0 \pmod 5$.
\end{itemize}
Then, the only possible solution of equation \eqref{eq1} is $n=2.$ 
\end{thm}
This theorem was claimed by Yuan and Zhang in \cite{X-Z} but their proof was incomplete. They missed the solution $(23, 78, 3)$ of  Lemma 2 in \cite{X-Z}, which corresponds to $D=7\cdot11\cdot79=6083$ in the paper.\\

Recently, Ishii \cite{I} gave a necessary and sufficient condition for equation \eqref{eq1} with the conditions $a \equiv 5 \pmod 6 ~\mbox{and}~ b \equiv 0 \pmod 3 $ to have positive integer solutions.

Let $d$ be a positive integer which is not a square. Then, by the theory of Pell equations, one knows the equation 
$$u^2-dv^2=1,~ u,v \in \mathbb{N}^+$$ 
has infinitely many positive integer solutions $(u,v)$ and all  of them are given by  
$$u_n+v_n\sqrt{d} = (u_1+v_1\sqrt{d})^n,$$ 
for some positive integer $n$, where $(u_1,v_1)$ is the smallest positive solution. This theory will help us in the statement of the next result.

\begin{thm} \label{thm2}
Suppose that $a$ is even and $b$ is a prime, $b \equiv 3 \pmod 8$. Then the equation 
$$(a^n-1)(b^n-1)=x^2$$  
has no solution in positive integers $(n,x).$  
\end{thm}
This result generalizes the main result of Szalay \cite{L}.
We organize this paper as follows. To prove the above results, we need some results on divisibility properties of the solutions of a Pell equation and some known results. These are presented in Section \ref{sec2}. The proof of Theorem \ref{thm1} is done in Section \ref{sec3}. We prove Theorem \ref{thm2} in Section \ref{sec4}.

\section{Preliminaries}\label{sec2}
In this section, we recall some results that are very useful for the proofs of our main results. The following result is well-known and one can refer to Lan and Szalay (see Lemma 1 of \cite{L-L}).
\begin{lem}\label{lem1}
Let $d$ be a positive which is not a square.
\begin{enumerate}
\item If $k$ is even, then each prime factor $p$ of $u_k$ satisfies $p \equiv \pm 1 \pmod 8.$
\item If $k$ is odd, then $u_1 \mid u_k ~\mbox{and}~ u_k/u_1$ is odd.
\item If $q \in \{2,3,5\}$, then $q \mid u_k$ implies $q \mid u_1.$
\end{enumerate}
\end{lem}

The following lemma can be deduced from \cite{B} and \cite{R}.
\begin{lem}\label{lem2}
Let $p>3$ be a prime. Then, the equation 
$$x^p=2y^2-1, \quad x, y \in \NN$$
 has only the solution $(x,y)=(1,1)$ in positive integers and the equation 
 $$x^3=2y^2-1, \quad x, y \in \NN$$
  has only the solutions $(x,y)=(1,1),(23,78)$ in positive integers.
\end{lem}

The last result to recall is Lemma 2.1 of \cite{G}.
\begin{lem}\label{lem3}
For a fixed $d$, if $2\mid u_r$ and $2\nmid u_s$, then  $2\nmid r$ and $2\mid s$.
\end{lem}
 
 \section{Proof of Theorem  \ref{thm1}}\label{sec3}

We prove only the first part of the statement, the proofs of the other parts are similar and left to the reader.

Let $a \equiv 2 \pmod 3$ and $b \equiv 0 \pmod 3$ and suppose that $(n,x)$ is a solution to equation \eqref{eq1}. Put $D=(a^n-1, b^n-1)$. From this equation, we have $$a^n-1=Dy^2, ~b^n-1=Dz^2, ~x=Dyz, ~ D,y,z \in \mathbb{N}.$$
Since $3\mid b,$ from $b^n-1=Dz^2$, it follows that 
$$D \equiv -1 \pmod 3  \quad \mbox{and} \quad 3\nmid z.$$ 
Now, we consider two cases according to weather $3$ divides $y$ or not. 

{\bf Case 1:} Suppose that $3\nmid y$. Then $y^2 \equiv 1 \pmod 3$ and we get 
$$a^n \equiv Dy^2+1 \equiv D+1 \equiv 0 \pmod 3.$$ 
This contradicts the fact that $a \equiv 2 \pmod 3$. 

{\bf Case 2:} Assume now that $3 \mid y.$ Since $a \equiv 2 \pmod 3$, from $a^n-1=Dy^2,$ we obtain 
$$2^n \equiv a^n \equiv Dy^2+1 \equiv 1  \pmod 3.$$ 
We deduce that $n$ is even. Put $n=2m$. Therefore, $D$ cannot be a square and the corresponding Pell equation $u^2-Dv^2=1$ has two solutions 
$$(u, v)=(a^m,y),(b^m,z).$$ 
Since $a\neq b,$ there exist distinct positive integers $r$ and $s$ such that 
$$(a^m,y)=(u_r,v_r) \quad \mbox{and} \quad (b^m,z)=(u_s,v_s),$$ 
where $(u_1,v_1)$ is the fundamental solution of this Pell equation. From Lemma \ref{lem1} and since $3\mid b,$ one can see that $2\nmid s$ and $3 \mid u_1.$ On the other hand, $a \equiv 2 \pmod 3$ and $3 \mid u_1$ imply that $2\mid r$  by Lemma \ref{lem1} (2). Put $r=2t$. Then we get 
$$a^m=u_{2t}=2u_t^2-1.$$ 
Here again, we consider two cases according to the parity of $m$. 

{\bf Case $2.1$:} Suppose that $2 \mid m$. Then $4\mid n$ and RESULT 2 in \cite{C} implies that $(a,b)=(13, 239)$. This contradicts the fact that $3\mid b.$ 

{\bf Case $2.2$:} We assume that $2\nmid m$. If $m>3,$  then since $m$ is odd, either $3 \mid m$ or $p \mid m$ for some $p>3$ prime. If $p \mid m$ for some $p>3$ prime, then Lemma \ref{lem2} shows that we have a contradiction since $a>1$.  If $3 \mid m$, then $a^{m/3}=23$ (by Lemma 4 for the prime exponent $3$) and $m/3 >1$ gives a contradiction since $23$ is not a perfect power of exponent greather than $1$ of some other integer. So $m/3 = 1$  and we obtain $a=23$ and $u_t=78$. As $(u_t,v_t)$ is a solution of Pell equation $x^2-Dy^2=1,$ it follows that 
$$Dv_t^2=u_t^2-1=78^2-1=7\times 11\times 79,$$ 
which shows that $D=7\times 11\times 79=6083$ but also that $v_t=t=1$ and $u_1=78.$ In particular, $2\| u_1.$ Since $s$ is odd and $u_s/u_1$ is odd, it follows that $2\|u_s,$ so $u_s$ cannot be the cube of $b$. So $m=1$, i.e., $n=2$. This completes the proof of Theorem  \ref{thm1}.\\


\section{Proof of Theorem  \ref{thm2}}\label{sec4}

Before starting the proof, let us remind the reader that a primitive divisor of $v_n$ is a prime factor $p$ of $v_n$ which does not divide $v_m$ for any $m<n$. Carmichael  \cite{Ca} proved that a primitive divisor exists whenever $n>6$. His result will be used below in our proof.

Suppose that $(a^n-1)(b^n-1)=x^2$ has solutions in positive integers $n,x$ when $a$ is even and $b$ is a prime such that $b\equiv 3 \pmod 8$.  Then we have 
$$a^n-1=Dy^2 \quad \mbox{and} \quad b^n-1=Dz^2,$$ 
where $D=(a^n-1, b^n-1)$. $D$ can be written into the form $D=dw^2,$ where $d$ is a square-free integer. 
If $n$ is odd, then $2||(b^n-1) = (b- 1)(b^{n-1}+\cdots +1)$,
showing that $2$ must divide $D$, which is a contradiction with the fact that $a$ is even.
So, $n$ is even. Thus, $n = 2m$ and $b^m = u_s$. Since $b\not\equiv \pm 1\pmod 8$, Lemma 3 tells us that $s$ is odd.
Since $u_1 | u_s$, it follows that $u_1 = b^{m_1}$, for some
$m_1\leq m$. We next show that $s = 1$. Assume that $s > 1$. It then follows that
$v_{2s} = 2v_su_s = 2v_sb^m$ has no primitive divisors 
since then all prime factors of
$v_{2s}$ are either in $\{b,2\}$ (which already divide $v_2 = 2u_1v_1$), or they  divide $v_s$.
By Carmichael's Primitive Divisor Theorem in \cite{Ca}, $2s\leq 12$ and since $s$ is odd, we
get $s\in \{3, 5\}$. Since $u_3 = 4u_1^3-3u_1$ and $u_5 = 16u_1^5-20u_1^3+5u_1$ and both $u_1$
and $u_3$ (or $u_1$ and $u_5$, respectively) are powers of the same prime $b$, we get
that $b\in \{ 3 ,5\}$. Now a result in \cite{LU}, shows that $s=1$, which is a contradiction to the assumption $s> 1.$ Thus, $u_1=b^m.$ Since $a$ is even, it follows that $a^m=u_r$ with odd $r$ (by Lemma 5). Since $u_r/u_1$ is odd (by Lemma 3), it follows that $2|u_1$. This contradicts the fact that $u_1=b^m$ and completes the proof of Theorem  \ref{thm2}. 
 
\section*{Acknowledgements} 
The authors thank the referee for a careful reading of the manuscript and for comments which improved its quality.  They are also particularly grateful for his/her suggestion for the proof of Theorem  \ref{thm2}. The second author is partially supported by Purdue University Northwest. The third author was supported by NSF of China (No. 11601476) and the Guangdong Provincial Natural Science Foundation (No. 2016A030313013 ).

 \bibliographystyle{plain}

\end{document}